\documentclass[12pt]{amsart}
 \usepackage[margin=1in]{geometry}
 \usepackage{amsmath,amsfonts,amsthm,amssymb,bbm}
 \usepackage{graphicx,color,dsfont}
 \usepackage{enumitem}
 \usepackage{fourier}
  \usepackage{hyperref}

\begin{document}

\newtheorem{theorem}{Theorem}
\newtheorem{lemma}{Lemma}
\newtheorem{proposition}{Proposition}
\newtheorem{rmk}{Remark}
\newtheorem{example}{Example}
\newtheorem{exercise}{Exercise}
\newtheorem{definition}{Definition}
\newtheorem{corollary}{Corollary}
\newtheorem{notation}{Notation}
\newtheorem{claim}{Claim}

\newtheorem{dif}{Definition}

 \newtheorem{thm}{Theorem}[section]
 \newtheorem{cor}[thm]{Corollary}
 \newtheorem{lem}[thm]{Lemma}
 \newtheorem{prop}[thm]{Proposition}
 \theoremstyle{definition}
 \newtheorem{defn}[thm]{Definition}
 \theoremstyle{remark}
 \newtheorem{rem}[thm]{Remark}
 \newtheorem*{ex}{Example}
 \numberwithin{equation}{section}

\newcommand{\vertiii}[1]{{\left\vert\kern-0.25ex\left\vert\kern-0.25ex\left\vert #1
    \right\vert\kern-0.25ex\right\vert\kern-0.25ex\right\vert}}

\newcommand{\R}{{\mathbb R}}
\newcommand{\C}{{\mathbb C}}
\newcommand{\U}{{\mathcal U}}
\newcommand{\cG}{W}
\newcommand{\tcG}{\tilde{\cG}}

\newcommand{\norm}[1]{\left\|#1\right\|}
\renewcommand{\(}{\left(}
\renewcommand{\)}{\right)}
\renewcommand{\[}{\left[}
\renewcommand{\]}{\right]}
\newcommand{\f}[2]{\frac{#1}{#2}}
\newcommand{\im}{i}
\newcommand{\cl}{{\mathcal L}}
\newcommand{\ck}{{\mathcal K}}

\newcommand{\al}{\alpha}
\newcommand{\be}{\beta}
\newcommand{\wh}[1]{\widehat{#1}}
\newcommand{\ga}{\gamma}
\newcommand{\Ga}{\Gamma}
\newcommand{\de}{\delta}
\newcommand{\ben}{\beta_n}
\newcommand{\De}{\Delta}
\newcommand{\ve}{\varepsilon}
\newcommand{\ze}{\zeta}
\newcommand{\Th}{\Theta}
\newcommand{\ka}{\kappa}
\newcommand{\la}{\lambda}
\newcommand{\laj}{\lambda_j}
\newcommand{\lak}{\lambda_k}
\newcommand{\La}{\Lambda}
\newcommand{\si}{\sigma}
\newcommand{\Si}{\Sigma}
\newcommand{\vp}{\varphi}
\newcommand{\om}{\omega}
\newcommand{\Om}{\Omega}
\newcommand{\ra}{\rightarrow}

\newcommand{\ro}{{\mathbf R}}
\newcommand{\rn}{{\mathbf R}^n}
\newcommand{\rd}{{\mathbf R}^d}
\newcommand{\rmm}{{\mathbf R}^m}
\newcommand{\rone}{\mathbf R}
\newcommand{\rtwo}{\mathbf R^2}
\newcommand{\rthree}{\mathbf R^3}
\newcommand{\rfour}{\mathbf R^4}
\newcommand{\ronen}{{\mathbf R}^{n+1}}
\newcommand{\ku}{\mathbf u}
\newcommand{\kw}{\mathbf w}
\newcommand{\kf}{\mathbf f}
\newcommand{\kz}{\mathbf z}

\newcommand{\N}{\mathbf N}

\newcommand{\tn}{\mathbf T^n}
\newcommand{\tone}{\mathbf T^1}
\newcommand{\ttwo}{\mathbf T^2}
\newcommand{\tthree}{\mathbf T^3}
\newcommand{\tfour}{\mathbf T^4}

\newcommand{\zn}{\mathbf Z^n}
\newcommand{\zp}{\mathbf Z^+}
\newcommand{\zone}{\mathbf Z^1}
\newcommand{\zz}{\mathbf Z}
\newcommand{\ztwo}{\mathbf Z^2}
\newcommand{\zthree}{\mathbf Z^3}
\newcommand{\zfour}{\mathbf Z^4}

\newcommand{\hn}{\mathbf H^n}
\newcommand{\hone}{\mathbf H^1}
\newcommand{\htwo}{\mathbf H^2}
\newcommand{\hthree}{\mathbf H^3}
\newcommand{\hfour}{\mathbf H^4}

\newcommand{\cone}{\mathbf C^1}
\newcommand{\ctwo}{\mathbf C^2}
\newcommand{\cthree}{\mathbf C^3}
\newcommand{\cfour}{\mathbf C^4}
\newcommand{\dpr}[2]{\langle #1,#2 \rangle}
\newcommand{\sech}{\textup{sech}}

\newcommand{\sn}{\mathbf S^{n-1}}
\newcommand{\sone}{\mathbf S^1}
\newcommand{\stwo}{\mathbf S^2}
\newcommand{\sthree}{\mathbf S^3}
\newcommand{\sfour}{\mathbf S^4}

\newcommand{\lp}{L^{p}}
\newcommand{\lppr}{L^{p'}}
\newcommand{\lqq}{L^{q}}
\newcommand{\lr}{L^{r}}
\newcommand{\echi}{(1-\chi(x/M))}
\newcommand{\chip}{\chi'(x/M)}

\newcommand{\wlp}{L^{p,\infty}}
\newcommand{\wlq}{L^{q,\infty}}
\newcommand{\wlr}{L^{r,\infty}}
\newcommand{\wlo}{L^{1,\infty}}

\newcommand{\lprn}{L^{p}(\rn)}
\newcommand{\lptn}{L^{p}(\tn)}
\newcommand{\lpzn}{L^{p}(\zn)}
\newcommand{\lpcn}{L^{p}(\cn)}
\newcommand{\lphn}{L^{p}(\cn)}

\newcommand{\lprone}{L^{p}(\rone)}
\newcommand{\lptone}{L^{p}(\tone)}
\newcommand{\lpzone}{L^{p}(\zone)}
\newcommand{\lpcone}{L^{p}(\cone)}
\newcommand{\lphone}{L^{p}(\hone)}

\newcommand{\lqrn}{L^{q}(\rn)}
\newcommand{\lqtn}{L^{q}(\tn)}
\newcommand{\lqzn}{L^{q}(\zn)}
\newcommand{\lqcn}{L^{q}(\cn)}
\newcommand{\lqhn}{L^{q}(\hn)}

\newcommand{\lo}{L^{1}}
\newcommand{\lt}{L^{2}}
\newcommand{\li}{L^{\infty}}
\newcommand{\beqn}{\begin{eqnarray*}}
\newcommand{\eeqn}{\end{eqnarray*}}
\newcommand{\pplus}{P_{Ker[\cl_+]^\perp}}

\newcommand{\co}{C^{1}}
\newcommand{\ci}{C^{\infty}}
\newcommand{\coi}{C_0^{\infty}}

\newcommand{\ca}{\mathcal A}
\newcommand{\cs}{\mathcal S}
\newcommand{\cm}{\mathcal M}
\newcommand{\cf}{\mathcal F}
\newcommand{\cb}{\mathcal B}
\newcommand{\ce}{\mathcal E}
\newcommand{\cd}{\mathcal D}
\newcommand{\cn}{\mathcal N}
\newcommand{\cz}{\mathcal Z}
\newcommand{\crr}{\mathbf R}
\newcommand{\cc}{\mathbf C}
\newcommand{\ch}{\mathcal H}
\newcommand{\cq}{\mathcal Q}
\newcommand{\cp}{\mathcal P}
\newcommand{\cx}{\mathcal X}
\newcommand{\eps}{\epsilon}

\newcommand{\pv}{\textup{p.v.}\,}
\newcommand{\loc}{\textup{loc}}
\newcommand{\intl}{\int\limits}
\newcommand{\iintl}{\iint\limits}
\newcommand{\dint}{\displaystyle\int}
\newcommand{\diint}{\displaystyle\iint}
\newcommand{\dintl}{\displaystyle\intl}
\newcommand{\diintl}{\displaystyle\iintl}
\newcommand{\liml}{\lim\limits}
\newcommand{\suml}{\sum\limits}
\newcommand{\ltwo}{L^{2}}
\newcommand{\supl}{\sup\limits}
\newcommand{\df}{\displaystyle\frac}
\newcommand{\p}{\partial}
\newcommand{\Ar}{\textup{Arg}}
\newcommand{\abssigk}{\widehat{|\si_k|}}
\newcommand{\ed}{(1-\p_x^2)^{-1}}
\newcommand{\tT}{\tilde{T}}
\newcommand{\tV}{\tilde{V}}
\newcommand{\wt}{\widetilde}
\newcommand{\Qvi}{Q_{\nu,i}}
\newcommand{\sjv}{a_{j,\nu}}
\newcommand{\sj}{a_j}
\newcommand{\pvs}{P_\nu^s}
\newcommand{\pva}{P_1^s}
\newcommand{\cjk}{c_{j,k}^{m,s}}
\newcommand{\Bjsnu}{B_{j-s,\nu}}
\newcommand{\Bjs}{B_{j-s}}
\newcommand{\Ly}{L_i^y}
\newcommand{\dd}[1]{\f{\partial}{\partial #1}}
\newcommand{\czz}{Calder\'on-Zygmund}
\newcommand{\chh}{\mathcal H}

\newcommand{\lbl}{\label}
\newcommand{\beq}{\begin{equation}}
\newcommand{\eeq}{\end{equation}}
\newcommand{\beqna}{\begin{eqnarray*}}
\newcommand{\eeqna}{\end{eqnarray*}}
\newcommand{\bp}{\begin{proof}}
\newcommand{\ep}{\end{proof}}
\newcommand{\bprop}{\begin{proposition}}
\newcommand{\eprop}{\end{proposition}}
\newcommand{\bt}{\begin{theorem}}
\newcommand{\et}{\end{theorem}}
\newcommand{\bex}{\begin{Example}}
\newcommand{\eex}{\end{Example}}
\newcommand{\bc}{\begin{corollary}}
\newcommand{\ec}{\end{corollary}}
\newcommand{\bcl}{\begin{claim}}
\newcommand{\ecl}{\end{claim}}
\newcommand{\bl}{\begin{lemma}}
\newcommand{\el}{\end{lemma}}
\newcommand{\dea}{(-\De)^\be}
\newcommand{\naa}{|\nabla|^\be}
\newcommand{\cj}{{\mathcal J}}

\title[Asymptotic attraction toward fronts of  the generalized Burgers equations]
{Asymptotic  attraction  with algebraic rates toward  fronts of  dispersive-diffusive  Burgers equations}

 \author[Milena Stanislavova]{\sc Milena Stanislavova}
 \address{Milena Stanislavova,   Department of Mathematics, University of Alabama - Birmingham,
 	1402 10th Avenue South
 	Birmingham AL 35294, USA
 }
 \email{mstanisl@uab.edu}
 \author[Atanas G. Stefanov]{\sc Atanas G. Stefanov}
 \address{Atanas G. Stefanov \\  Department of Mathematics, University of Alabama - Birmingham,
 	1402 10th Avenue South
 	Birmingham AL 35294, USA.
 }
 \email{stefanov@uab.edu}

\thanks{ Milena Stanislavova is partially supported by the National Science Foundation,  under award \# 2210867.  Atanas Stefanov   is partially  supported by  NSF-DMS under grant  \# 2204788.}

\subjclass[2010]{Primary  }


\date{\today}
 
\begin{abstract}
Burgers equation is a classic model, which arises in numerous applications. At  its very core it is a simple conservation law, which serves as a toy model for various dynamics phenomena. In particular, it supports explicit heteroclinic solutions, both fronts and backs.  Their stability has been studied in details. There has been substantial interest in considering dispersive and/or diffusive modifications, which present novel dynamical paradigms in such simple setting. More specificaly,  the KdV-Burgers model has been showed to support unique fronts (not all of them monotone!) with fixed values at $\pm \infty$. Many articles,  among which \cite{Pego}, \cite{NS1}, \cite{NS2}, have studied the question of stability of monotone (or close to monotone) fronts.
 
 In a breakthrough paper, \cite{BBHY}, the authors have extended these results in several different directions. They have considered a wider range of models. The fronts do not need to be monotone, but are subject of a spectral condition instead.  Most importantly the method allows for large perturbations, as long as the heteroclinic conditions at $\pm \infty$ are met. That is, there is asymptotic attraction to the said fronts or equivalently the limit set consist of one point. 
 
The purpose of this paper is to extend the results of  \cite{BBHY} by providing explicit algebraic rates of convergence as $t\to \infty$. We bootstrap these results from the results in \cite{BBHY} using additional energy estimates for two important examples namely KdV-Burgers and the fractional Burgers problem. These rates are likely not optimal, but we conjecture that they are algebraic nonetheless. 
 
\end{abstract}

\maketitle
\section{ Introduction}
Burgers equation 
\begin{equation}
	\label{b:10} 
	u_t - u_{xx} +u u_x=0, u(0,x)=u_0(x)
\end{equation}
where $u:\rone_+\times \rone \to \rone$   is a particularly simple example of a conservation law in one spatial dimension. As a simple model, which resembles (both in its form as well as its behavior) many distinguished equations in mathematical physics, it is  an ubiquitous object in the modern dynamical systems theory. The Cauchy problem, along with many other properties, has been explored in the literature. In particular, and specifically for this model, the Cole-Hopf transformation has played an important role as it transforms \eqref{b:10} into a linear problem.  More specifically, the  assignment $u=-2\f{w_x}{w}$ reduces \eqref{b:10} to the linear heat equation for $w$, namely $$w_t-w_{xx}=0, w(0,x)=e^{-\f{1}{2}\int_0^x u_0(y) dy}$$ In this way, one can write explicit solutions and derive particular properties of the conservation laws, which are not transparent in the usual general implicit formulations. 

 Along with the Burgers equation, related models have been analyzed. The KdV-Burgers model
  \begin{equation}
 	\label{b:12} 
 u_t-u_{xx}+u u_x= \nu  u_{xxx}, 
\end{equation}
   which describes, among other things, bores in viscous fluids, \cite{SG} and weak plasma shocks propagating transversely to a magnetic field. Many similar models were analyzed recently, among them the following general model, \cite{BBHY},
   \begin{equation}
   	\label{b:14} 
   	u_t - u_{xx}+u u_x = \cl u. 
   \end{equation}
Here $\cl$ is a multiplier operator, given by  $\widehat{\cl f}(k) = l(k) \hat{f}(k)$. Now, me must  impose  assumptions on   $\cl$. Namely, it is required that $\cl$ maps  real-valued functions into real-valued functions, and also its quadratic form is non-negative. That is, it satisfies 
\begin{eqnarray}
	\label{a:14} 
& & 	\dpr{\cl f}{f} =\Re \int_{-\infty}^\infty l(k) |\hat{f}(k)|^2 dk\leq 0 \\
\label{210} 
	& & \cl [1]=0
\end{eqnarray}
The condition \eqref{210} can be alternatively stated as $l[0]=0$ and $l$ is continuous at zero, but the form \eqref{210} is more suitable in this context. 
Next, one might argue that the condition \eqref{a:14} is akin to the dispersivity of the system (especially when one has that $\dpr{\cl f}{f}=0$), but as the main effects on the dynamics shall remain dissipative in nature, we will refrain from labeling  it as such.  
Note that \eqref{b:14}, incorporates \eqref{b:12}, but includes much more. In fact,  we list several physically relevant  cases: the KdV-Burgers, $\cl=\nu \p_{xxx}$,  the Benjamin-Ono-Burgers,  $\cl=\p_x |\p_x|$, the Hilbert-Burgers case $\cl= \ch=\p_x|\p_x|^{-1}$, the fractional Burgers equation\footnote{For the definition of the fractional derivatives please consult Section \ref{sec:2} below}  with $\cl=-\sum_{j=1}^N a_j (-\p_{xx})^{\al_j}, \ \ \ 0<\al_1<\al_2<\ldots <\al_N<1, a_j>0$ etc. All of these could be of interest.

 Understanding of the long time dynamics of such models is clearly important, and as it turns out a challenging problem.  As is usual in conservation laws, the standard background is starting with heteroclinic data with enough decay towards  the limiting values at $\pm \infty$. To this end, assume that the initial data $u_0$ satisfies 
 \begin{equation}
 	\label{200} 
 	\lim_{x\to \pm \infty} u_0(x)=u_\pm
 \end{equation}
 where $u_->u_+$ are constants.
 \subsection{Galilean transformation: reductions and existence of fronts}
 
 Let us record two elementary properties of the equation \eqref{b:14}. For the Burgers equation, i.e. the model \eqref{b:14} with $\cl=0$, we have true  invariances with respect to the Galilean transformation $$u(t,x)=v(t,x- ct)+c,$$ as well as the scaling transformation $$u\to u_\la:=\la u(\la^2 t, \la x), \la>0$$ 
 These two transformations,  in general,  do not commute with the flow of \eqref{b:14}, but they transform it into a model of the same type, i.e. it still satisfies \eqref{a:14}, \eqref{210}.   Since our goal is to understand the long term dynamics of heteroclinic initial data, we can use these transformations to our advantage as follows: given $u_{\pm}, u_->u_+$, introduce $c=\f{u_- + u_+}{u_--u_+}, \la=\f{u_- - u_+}{2}>0$ and a function $v$, so that 
 $$
 u(t,x)=\la U(\la^2 t, \la x - c\la^2 t) + c \la 
 $$
  Then, the function $U$ will satisfy the equation\footnote{Here note that the condition \eqref{210} is instrumental in establishing \eqref{230}, as $\cl$ annihilates the constants.}
  \begin{equation}
  	\label{230}
  	U_t - U_{xx}+U U_x=\cl_1  U, 
  \end{equation}
 where $\widehat{\cl_1 g}(k)=\la^{-2} l(\la k) \hat{g}(k)$, $\lim_{x\to \mp \infty} U(x) =\pm 1$.  As we shall see, $\cl_1$ also satisfies the relevant assumptions, \eqref{a:14}, \eqref{210}, so it will suffice to work with the model, where $U_{\mp}=\pm 1$. Naturally, and by the choice of the transformations, this corresponds to a case where $c=0$ as well. Thus, we restrict our attention to the steady state solutions of \eqref{b:14} with end values $\mp 1$, 
  which in this case are functions $\phi:\rone\to \rone$,  
  \begin{equation}
  	\label{p:10} 
  	-\phi''+\phi \phi'=\cl \phi; \ \ \lim_{x\to \mp \infty} \phi(x)=\pm 1. 
  \end{equation}
  The existence theory for such a problem is of course highly non-trivial. In fact, as in previous papers, we shall make assumptions to this effect. We mention in this regard that for the KdV-Burger's model, the existence of such monotone $\phi$ was proved in \cite{BS1, BVS}.  It is worth adding that these fronts are monotone if and only if $|\nu|\leq \f{1}{4}$, \cite{BVS}.  In addition, such solutions are unique (if one specifies the values $u_{\pm}$), and one can write the explicit form of the fronts. For example, in the case\footnote{and all the other cases can be recovered from the Galilean invariance} $\nu=-\f{1}{10}$, the formula is as follows 
  $$
  u(t,x)=\f{6}{5}\left[\sech^2\left(x+\f{12}{5}t\right)-
  2\tanh\left(x+\f{12}{5}t\right)-2\right].
  $$
   see \cite{JM}.

  \subsection{Asymptotic attraction toward heteroclinic fronts} 
We now review the literature regarding earlier asymtotic stability results. R. Pego, \cite{Pego} has shown that the {\it monotonic} fronts are orbitally asymtotically stable. More specifically, he showed that  any (small) perturbation of the front converges to some translate of the front. In later works, Naumkin-Shishmarev, \cite{NS1}, \cite{NS2} showed that the KdV-Burgers fronts are asymptotically stable, even for some  slightly non-monotonic fronts, i.e. $\nu: 0<\nu-\f{1}{4}<<1$.  

Next, we introduce the results of Barker-Bronski-Hur-Yang, \cite{BBHY}, which represents a real breakthrough in the study of this problem. Specifically, the authors improved the aforementioned results in at least three important directions. First, they consider all (non-necessarily monotonic) fronts, subject to a natural and specific spectral condition\footnote{which appears to be necesary}, see Theorem \ref{theo:BBHY} below. Second, they extend the range of applicability to  a wide range of examples, driven by possibly a mixture of dissipative and/or dispersive linear operators, subject to \eqref{a:14}. Third, and most impressively, they allow large data perturbations for asymptotic attraction type of result viz. a. viz. asymptotic stability. That is, as long as the initial data approaches the fixed homoclinics with good decay properties, and the fronts exist (with the prescribed spectral condition), the solution will stay close to translates (defined dynamically) of the front, all in $L^2$ norm. 
 	 
 	 We state the precise result below. 
 	 \begin{theorem}(Barker-Bronski-Hur-Yang, \cite{BBHY}) 
 	 	\label{theo:BBHY} 
 	 	
 	 	Suppose that $\cl$ satisfies the dispersivity property \eqref{a:14}, while the front $\phi$ satisfies  the profile equation \eqref{p:10} and in addition. 
 	 		\begin{equation}
 	 		\label{30} 
 	 		\lim_{x\to \mp \infty} \phi(x)=\pm 1, 	\phi', \phi''\in L^2(\rone), \int_{\rone} (1+|x|) |\phi'(x)| dx<\infty.
 	 	\end{equation}
  	Finally, assume that for some $\eps\in (0,1)$, the Schr\"odinger operator
  	\begin{equation}
  		\label{32} 
  		\ch_\eps:=-(1-\eps) \p_{x}^2+\f{1}{2} \phi' \ \ \textup{has exactly one negative eigenvalue}.
  	\end{equation} 
  Assume that $u_0$ is a heteroclinic data, so that 
  	$	\lim_{x\to \mp \infty} u_0(x)=\pm 1$ and more precisely, 
  	$u_0(x)-\phi\in L^2(\rone)$. 
  	Then, there exists a continuous function $x_0(t): \rone_+\to \rone$, so that the solution to \eqref{b:14} can be written in the form 
  	\begin{equation}
  		\label{n:20}
  		u(t, x)=\phi(x-x_0(t))+v(t, x-x_0(t)), 
  	\end{equation}
  	 so that $t\to \|v(t, \cdot)\|_{L^2(\rone)}$ is a monotonically decreasing function of time  and 
  	 \begin{equation}
  	 	\label{35} 
  	 	\lim_{t\to \infty} \|v(t, \cdot)\|_{L^p(\rone)}=0, p\in (2, \infty]. 
  	 \end{equation}
 	 \end{theorem}
  Let us record an energy estimate, which came up in the proof of Theorem \ref{theo:BBHY}. More precisely,  and under the assumptions of Theorem \ref{theo:BBHY}, and specifically \eqref{32}, we have 
  \begin{equation}
  	\label{50} 
  	\p_t \|v\|_{L^2}^2\leq -C \|v'\|_{L^2}^2. 
  \end{equation}
We shall make extensive use of \eqref{50} in the sequel, but note right away that the monotonicity of the mapping $t\to \|v(t, \cdot)\|_{L^2(\rone)}$ is a direct consequence of it. 

  Next, a few remarks are in order. 
  \begin{enumerate}
  	\item Note that the results in Theorem \ref{theo:BBHY} do not require smallness of the perturbation $u-\phi$. In other words, this is a strong asymptotic attraction result rather than a mere asymptotic stability result. 
  	\item Related to the previous point, this implies a rather strong uniqueness (up to the usual translation) theorem for  such fronts $\phi$. Namely, fronts $\phi$ with the properties \eqref{30}, \eqref{32}, if any,  are unique. The uniqueness result  was certainly explored in the literature in certain contexts (see for example \cite{JM} for the KdV-Burgers case), but not in the generality of the model \eqref{b:14}. 
  	\item  In the statement of Theorem \ref{theo:BBHY}, the particular values of the limits $\phi_{\pm}=\lim_{x\to \pm\infty} \phi(x)$, are set to be  $\phi_{\pm}=\mp 1$  only out of convenience. Assuming that solutions exists for some other pair , subject to $\phi_{\pm}: \phi_->\phi_+$, one can retrace the proof in \cite{BBHY} to the same results, as stated in Theorem \ref{theo:BBHY}. 
  	\item Note that since $\int_{-\infty}^{+\infty}\phi'(x) dx=-2$, the Schr\"odinger operators $\ch_\eps$ do have negative spectrum, which consists of eigenvalues. The condition \eqref{32} requires that the spectrum is minimal, in the sense that it consists of only one   eigenvalue. 
  	\item Related to the previous item, we note that the spectral condition \eqref{32} has been {\it rigorously} verified in \cite{BBHY},  for the (unique) monotone KdV-Burger's fronts, i.e. for the case $|\nu|<\f{1}{4}$. In the companion paper, \cite{BBHY2} the authors further explored the validity of \eqref{32} beyond the monotone cases. In fact, informal numerics show that \eqref{32}  holds all the way up to $|\nu|<4.1$, while rigorous numerics confirm this up to $|\nu|<3.9$. 
  \end{enumerate}

\subsection{Main result: pointwise decay for the KdV-Burgers equation}
Our first result concerns the standard KdV-Burgers model. We  bootstrap upon the results of \cite{BBHY} to obtain pointwise in time bounds for various Lebesgue norms. The precise statement is below. 
\begin{theorem}(Pointwise decay to the fronts for the KdV-Burgers model) 
	\label{theo:20} 
	For the KdV-Burgers model, i.e. $\cl=\nu \p_{xxx}$, under the assumptions of Theorem \ref{theo:BBHY} and  
	\begin{equation}
		\label{n:50}
			\int_{\rone} |u_0(x)-\phi(x)|^2 (1+|x|) dx<\infty, 
	\end{equation}
	the decomposition \eqref{n:20} still holds, with 
	\begin{equation}
		\label{118} 
		\|v(t, \cdot)\|_{L^2}\leq \f{C}{\sqrt{t}}. 
	\end{equation}
	In addition, for every $q\in (1,2)$, there is a decay $\lim_{t\to \infty} \|v(t, \cdot)\|_{L^q}=0$, and in fact for every $\de>0$, there exists $C=C_\de$, so that 
	\begin{eqnarray}
		\label{120} 
		\|v(t, \cdot)\|_{L^q}\leq C_\de t^{-(1-\f{1}{q})+\de}, 1<q<2, 
	\end{eqnarray}
	For $p>2$,  the bound is 
	\begin{equation}
		\label{124} 
		\|v(t, \cdot)\|_{L^p}\leq C  t^{-\f{1}{p}}, 2\leq  p<\infty.  
	\end{equation}
\end{theorem}
Some important remarks to be made herein:
\begin{enumerate}
	\item While we do require a more localized data,   see \eqref{n:50}, we still allow such perturbations to be large.
	\item The decay rate \eqref{118} is faster, compared to the linear heat equation, where we only have $\|v(t)\|_{L^2(\rone)}\sim t^{-\f{1}{4}}$. This  enhance decay rate can be thought of as  an anomalous dissipation of sorts,  caused by the non-trivial background around the front solution. 
\end{enumerate}

\subsection{Main result: Pointwise decay for the fractional Burgers equation}
For clarity, we consider the fractional Burgers equation to be the evolution problem
\begin{equation}
	\label{fB} 
	u_t - u_{xx}+\sum_{j=1}^N a_j (-\p_{xx})^{\al_j} u + u u'=0,
\end{equation}
subject to initial data $u(0,x)=u_0(x)$. Clearly, the evolution of \eqref{fB} preserves odd initial data. We have the following, more precise  than \eqref{35} result. 
\begin{theorem}
	\label{theo:10} 
	Suppose that $\phi$ is an {\bf odd} stationary front  of the fractional Burger equation 
	\begin{equation}
		\label{fb:10} 
		-\phi''+\sum_{j=1}^N a_j (-\p_{xx})^{\al_j} \phi +\phi\phi'=0.
	\end{equation}
where $0<\al_1<\al_2<\ldots <\al_N<1, a_j\geq 0$. In addition, $\phi$ satisfies \eqref{30}and the spectral  assumption as in Theorem \ref{theo:BBHY}. 
Suppose that $u_0=\phi+v_0$, where $v_0\in L^1\cap L^2(\rone)$ is an odd function. 
  Then, the solution $u$ of \eqref{b:14}, with initial data $u_0$ can be written as $$u(t,x)=\phi(x)+v(t,x),$$ where  for $t>2$, one has the bounds 
  \begin{eqnarray}
  	\label{90} 
  	\|v(t, \cdot)\|_{L^p} &\leq & C \|v_0\|_{L^2\cap L^1} 
  	\left(\f{\ln(t)}{t}\right)^{\f{1-\f{1}{p}}{2}}, 1<p<2\\
  	\label{92} 
  	\|v(t, \cdot)\|_{L^q} &\leq & C \left(\f{\ln(t)}{t}\right)^{\f{7}{24}-\f{1}{12 q}}, 2<q\leq \infty
  \end{eqnarray}
\end{theorem}
{\bf Remarks:} 
\begin{itemize}
	\item Note that as before, the perturbation $v_0$ may be large, as long as it is well-localized, i.e. $v_0\in L^1\cap L^2(\rone)$.  
	\item Related to the previous item, it is worth noting that our proof fails under an  assumption such as $v_0\in L^q\cap L^2, 1<q$. That is, we cannot establish the algebraic decay rates (for any $p\in (1, \infty]$!) in \eqref{90}, \eqref{92}, without assuming first that  $v_0\in L^1(\rone)$. At this point it is unclear whether this is a technical obstacle or one indeed needs to require  $L^1$ decay of the perturbation. 
	\item As before, Theorem \ref{theo:10} provides strong uniqueness property of the odd fronts $\phi$ satisfying \eqref{fb:10}. 
	\item It is an open question whether or not the algebraic decay rates stated in \eqref{90} are sharp. While the logarithmic corrections are most certainly not, one might wonder about the other bounds in \eqref{90}, \eqref{92}. 
\end{itemize}
We plan our paper as follows. In Section \ref{sec:2}, we introduce some basic notations, as well as the Sobolev embedding and various related interpolation inequalities. We introduce the notion of a diffusive semigroup, as this allows us for an effective energy estimates for the fractional Burgers evolution in $L^p, 1<p$ norms. In Section \ref{sec:3}, we consider the KdV-Burgers problem and we show enhanced energy estimates (for the weighted norms $t\to \left(\int_{\rone} |v(t,x)|^2 |x| dx\right)^{\f{1}{2}}$), which allows us to bootstrap the results of \cite{BBHY} to pointwise in time decay estimates for $\|v(t, \cdot)\|_{L^2}$. In Section \ref{sec:4}, we use the integration by parts techniques of \cite{CL} to derive energy estimates for $\|v(t, \cdot)\|_{L^p}, p>1$ and ultimately for $p=1$. We again use in a critical way the energy estimate \eqref{50} from \cite{BBHY}, and in conjunction with the {\it a priori} estimates for $\|v(t, \cdot)\|_{L^1}$, we conclude pointwise in time decay for the solution norms. 
\section{Preliminaries} 
\label{sec:2} 
We shall use standard notations for the Lebesgue spaces $L^p, 1\leq p\leq \infty$. More generally, the Sobolev spaces $W^{k,p}(\rone)$ for integer $k, 1<p<\infty$, are defined via the norms 
$$\|f\|_{W^{k,p}}=\|f\|_{L^p}+ \sum_{l=1}^k \|\p_x^l f\|_{L^p}$$
\subsection{Fourier transform, fractional Burger's and the Bernstein inequality}
We take the Fourier transform and its inverse in the form 
$$
\hat{f}(\xi)=\int_{\rone} f(x) e^{-2\pi i x\xi} dx, \ \ f(x)=\int_{\rone} \hat{f}(\xi) e^{2\pi i x\xi} d\xi.
$$
In this notations, we have $\widehat{f''}(\xi)=-4\pi^2 \xi^2 \hat{f}(\xi)$. More generally, we introduce the fractional operators $(-\p_{xx})^\al$ as acting on Schwartz functions via the symbol $(2\pi|\xi|)^{2\al}$. More precisely, 
$$
\widehat{(-\p_{xx})^\al f}(\xi)= (2\pi|\xi|)^{2\al} \hat{f}(\xi).
$$
In fact, this allows one to introduce a companion spaces to $W^{s,p}$, where $s>0, 1<p<\infty$ is not necessarily an integer. More precisely, these are given by the norm 
$$
\|f\|_{W^{s,p}}=\|f\|_{L^p}+\| (-\p_{xx})^{\f{s}{2}} f\|_{L^p}.
$$
These spaces are equivalent to the previously defined $W^{p,k}$ for integer $s$. The Sobolev inequality in one spatial dimension takes the simple form 
$$
\|u\|_{L^q(\rone)}\leq C \|u\|_{W^{s,p}} \ \ 1<p<q<\infty, \ \ s\geq \f{1}{p}-\f{1}{q}. 
$$
Next, we discuss the Bernstein's inequality, which is a relatively crude, but helpful form of Sobolev embedding. Consider a function $f: \rn\to \rn$ and a set of finite Lebegue measure $A$. Let $f_A: \widehat{f_A}(\xi)=\hat{f}(\xi) \chi_A(\xi)$,  where  $\chi_A(x)=\left\{\begin{array}{cc}
	1 & x\in A\\
	0 & x\notin A
\end{array}
\right..$
Then, for all $1\leq p\leq \infty$, 
\begin{equation}
	\label{75} 
	\|f_A\|_{L^q(\rn)}\leq |A|^{\f{1}{p}-\f{1}{q}} \|f\|_{L^p(\rn)}.
\end{equation}
An easy estimate, which follows  immediately from Plancherel's identity is 
$$
\|\p_x f_A\|_q\leq C \min_{\xi\in A}  |\xi| \|f_A\|_q, 1\leq q\leq \infty.
$$
We will use the rough Fourier cutoff projections
$$
f_{<\eps}: \widehat{f_{<\eps}}(\xi):=\hat{f}(\xi)\chi_{(-\eps, \eps)}(\xi), \ \ f_{>\eps}:=f- f_{<\eps}.
$$
We also need the Gagliardo-Nirenberg estimate (or  the log-convexity of the map $p\to \|f\|_p$).  More precisely, for every $1\leq p<q<r\leq \infty: \f{1}{q}=\f{\theta}{p}+\f{1-\theta}{r}$, 
$$
\|f\|_q\leq \|f\|_p^\theta \|f\|_r^{1-\theta}.
$$
Finally, there is the estimate, 
$
\|f\|_{L^\infty}^2 \leq C \|f'\|_{L^2} \|f\|_{L^2},
$
which will also be useful. 
\subsection{Diffusion semigroups} 
The notion of diffusive semigroups will play a role in our considerations. Shortly, we say that  $\{T(t)\}_{t\geq 0}$ is  a $C_0$ semigroup
on a Banach space, if $T(t)\in B(X)$, $T(t+s)=T(t)T(s), T(0)=Id$ and for every $f\in X: \lim_{t\to 0+} \|T(t) f - f\|_X=0$.  It is standard that one can always associate with such objects 
a generator $L$, a generally unbounded operator, which can be introduced on a domain $D(L)=\{f\in X: \lim_{t\to 0+} \frac{T(t) f - f}{t} \ \  \textup{exists}\}$ as 
$$
L f:=\lim_{t\to 0+} \frac{T(t) f - f}{t}, f\in D(L).
$$
Naturally,  $T(t) f\in D(L)$ solves the $X$ valued differential equation $x'=L x, x(0)=f$ for any $f\in X$, so we informally denote it by $T(t)=e^{t L}$. 
\begin{definition}
	\label{defi:10} We say that  a $C_0$ semigroup on $L^2(\rn)$, $T(t)$,   is a diffusive semigroup, if  $T(t), t\geq 0$ is  self-adjoint, and given by a  probability density  kernel, i.e. 
	$$
	T(t) f(x) = \int_{\rn} p_t(x-y) f(y) dy, f\in L^2(\rn), p_t(x)\geq 0, \int_{\rn} p_t(x) dx = 1.
	$$
\end{definition} 
We observe that the generator $L$ in such cases is necessarily a non-positive operator, since for all $t\geq 0$ and a test function $f$, we have 
\begin{eqnarray*}
	\dpr{T(t) f}{f} &=& \int  p_t(x-y) f(y) \bar{f}(x) dx dy  \leq  \\
& \leq & 	\left(\int p_t(x-y) |f (y)|^2 dx dy\right)^{\f{1}{2}} \left(\int p_t(x-y) |f (x)|^2 dx dy\right)^{\f{1}{2}}=\|f\|^2,
\end{eqnarray*}
whence $\dpr{L f}{f} \leq \limsup_{t\to 0+} \frac{\dpr{T(t) f - f}{f}}{t}\leq 0$. 
 
It is standard that the operator $L=\p_{xx}$ is a generator of such semigroup. More generally, the fractional Burgers operators $L=-\sum_{j=1}^N a_j (-\p_{xx})^{\al_j}, a_j\geq 0, 0<\al_1<\ldots \al_N\leq 1$ also generate diffusion semigroups\footnote{Here the assumption $\al_j\leq 1$ is crucial, since it is well-known that the kernels of $e^{-t (-\p_{xx})^\al}$ are sign-changing, as soon as $\al>1$}, see for example \cite{CL}. An important property of such semigroups is the following. 
\begin{proposition}(Theorem 3.2, \cite{CL})
\label{prop:16} 
Assume that  $L$ generates a diffusive  semigroup. Then, 
for all $1\leq p<\infty$, we have 
\begin{equation}
		\label{40}
  	\int_{\rn} f|f|^{p-2} L(f) dx\leq 0. 
\end{equation}
If, in addition, $ 2\leq p<\infty$, then 
\begin{equation}
	\label{39} 
	 	-p \int f|f|^{p-2} L(f) dx \geq \int |\sqrt{-L}(f|f|^{\f{p}{2}-1})|^2 dx
\end{equation}
\end{proposition}
In particular, we have that for each $0<\al\leq 1, 2\leq p<\infty$
\begin{equation}
	\label{42} 
		\int_{\rn} ((-\De)^\al f) f|f|^{p-2} dx \geq 0.
\end{equation}

\section{Attraction to fronts: the Burger-KdV model}
\label{sec:3}

We start by revisiting some of the results obtained in the proof of Theorem \ref{theo:BBHY}.   Indeed, in the course of the proof, the authors set $\dot{x}_0(t)=\ga \int_{-\infty}^\infty \phi'(x) v(x) dx$, with\footnote{But if one uses different values of $\phi_\pm$, the value of $\ga$ must change according to $\ga>\f{2}{\phi_- -\phi_+}$ } $\ga>1$.  This is done to ensure that the rank-one perturbation $\ch_\ga=-\p_{xx}+\f{1}{2}\phi'+\ga \dpr{\phi'}{\cdot} \phi'$ of the important linearized operator\footnote{Recall that $\ch$, and even slightly smaller operator $-(1-\eps)\p_{xx}+\f{1}{2}\phi'$ have a negative e-value}   $\ch=-\p_{xx}+\f{1}{2}\phi'$ has the property $\ch_\ga\geq 0$.  Two important features of the construction  are 
\begin{eqnarray}
	\label{129} 
& & t\to \|v(t, \cdot)\|_{L^2}^2 \ \ \textup{is monotonically decreasing}  \\  
	\label{130} 
& & 	\int_0^t |\dot{x}_0(s)|^2 ds+	\int_0^t  \|v'(s, \cdot)\|_{L^2}^2 ds \leq C \|v_0\|_{L^2}^2 
\end{eqnarray}
We make substantial use of these two facts in our subsequent analysis. Next, we write the equation for the residual $v: u(x,t)=\phi(x-x_0(t))+v(t, x-x_0(t))$. This takes the form 
\begin{equation}
	\label{140} 
	v_t - v_{xx} -\dot{x}_0 v_x-\dot{x}_0 \phi'+\phi' v+\phi v'+ v v'+\nu v_{xxx}=0.
\end{equation}
We setup an energy estimate for $v$ in the space $L^2_1(\rone)$. More precisely, we shall show that 
$$
\int_{-\infty}^\infty |v(t,x)|^2 |x| dx\leq C, t\geq 0.
$$
Of course, the global solution $v$ produced in \cite{BBHY} does not {\it a priori} live in the space $L^2_1(\rone)$, so one begins by only  setting formally the energy estimate for it. We will ignore the details of this justification, as this is standard by considering instead a quantity of the form $\int_{-\infty}^\infty |v(t,x)|^2 |x| \zeta(x/N)dx$ for $N>>1$ and some smooth cutoff function $\zeta$. With these caveats, we proceed to bound the $L^2_1(\rone)$ norm. Taking a dot product of \eqref{140} with $v |x|$, we have
\begin{eqnarray*}
	& & \f{1}{2} \p_t  \int v^2(t,x) |x| dx  -\int v'' v |x| dx   -\dot{x}_0 \int v' v |x| dx - 
	\dot{x}_0 \int \phi' v |x| dx+\int \phi' v v |x| dx +\\
	&+& \int \phi v' v |x| dx + \int v v' v |x| dx +\nu \int v''' v |x| dx =0.
\end{eqnarray*}
We now proceed to estimate various terms, by keeping in mind the relation \eqref{130}, which guides us about  the acceptable error terms. 
After integration by parts and algebraic manipulations, we have
\begin{eqnarray*}
	-\int v'' v |x| dx &= & \int (v')^2 |x| dx + \int_{-\infty}^\infty  v' v sgn(x) dx 
\geq -v^2(t,0)\geq -\|v(t, \cdot)\|_{L^\infty}^2,\\
  \left| \dot{x}_0 \int v' v |x| dx \right| &=& \left| \f{\dot{x}_0}{2}  \int v^2 sgn(x)   dx \right| \leq C |\dot{x}_0 (t)| \|v(t, \cdot)\|_{L^2}^2\leq  C |\dot{x}_0 (t)| \|v(t, \cdot)\|_{L^2}\\ 
  \left|	\dot{x}_0 \int \phi' v |x| dx\right|&\leq& C |\dot{x}_0 (t)|  \|v(t, \cdot)\|_{L^2}, \\
  \left|\int v v' v |x| dx\right|&=&\f{1}{3} \left|\int v^3 sgn(x)dx\right|
\leq C \|v(t, \cdot)\|_{L^3}^3\leq C \|v\|_{H^{\f{1}{6}}}^3\leq C \|v'(t, \cdot)\|_{L^2}^{\f{1}{2}} 
\|v(t, \cdot)\|_{L^2}^{\f{5}{2}}\\
&\leq & C \|v'(t, \cdot)\|_{L^2}^{\f{1}{2}} 
\|v(t, \cdot)\|_{L^2}^{\f{3}{2}}
\end{eqnarray*}
Also, 
$$
\int \phi' v^2 |x| dx\geq -\|v(t, \cdot)\|_{L^\infty}^2 \int |\phi'(x)|  |x| dx =-C \|v(t, \cdot)\|_{L^\infty}^2. 
$$
Here,  in addition to Sobolev embedding and Gagliardo-Nirenberg estimates, we   used the bound $\|v(t, \cdot)\|_{L^2}\leq C$.
Next, 
\begin{eqnarray*}
	& & 
\int \phi v' v |x| dx=-\f{1}{2}   \left( \int \phi'(x) v^2(x) |x| dx + \int \phi(x) v^2(t,x) sgn(x) dx \right)\geq \\
& & \geq  - C  \|v(t, \cdot)\|_{L^\infty}^2-\f{1}{2}\int \phi(x) v^2(t,x) sgn(x) dx 
\end{eqnarray*}
Now, note that  $-\phi(x) sgn(x)>0$ for all $x$ and also there exists $A$, so that $-\phi(x) sgn(x)>\f{1}{2}$ for $|x|>A$. This last statement follows from $\lim_{x\to \pm \infty} \phi(x) = \mp 1$. So, we further estimate
\begin{eqnarray*}
& & 	-\f{1}{2}\int \phi(x) v^2(t,x) sgn(x) dx> -\f{1}{2}\int_{|x|>A}  \phi(x) v^2(t,x) sgn(x) dx >  \f{1}{4} \int_{|x|>A} v^2(t,x) dx =  \\
& & = \f{1}{4} \int v^2(t,x) dx -\f{1}{4}  \int_{-A}^A v^2(t,x) dx \geq    \f{1}{4} \int v^2(t,x) dx - C  \|v(t, \cdot)\|_{L^\infty}^2, 
\end{eqnarray*}
Finally, for the term $\int v''' v |x| dx$, we have 
$$
\left|\int v''' v |x| dx \right|= \left| \f{1}{2} \int (v')^2 sgn(x) - \int v'' v sgn(x) dx \right| \leq 
2 |v(0)| |v'(0)| + C \|v'(t, \cdot)\|_{L^2}^2. 
$$
Putting it all together and estimating $-\|v(t, \cdot)\|_{L^\infty}^2 \geq - C\|v'(t, \cdot)\|_{L^2} \|v(t, \cdot)\|_{L^2}$, we obtain 
\begin{eqnarray*}
& & 	\p_t  \int v^2(t,x) |x| dx + \f{1}{4} \int v^2(t,x) dx \leq C( |\dot{x}_0 (t)| \|v(t, \cdot)\|_{L^2}+ |v(t,0)| |v'(t,0)| )+ \\
	& & + C( \|v'(t, \cdot)\|_{L^2}^2+\|v'(t, \cdot)\|_{L^2} \|v(t, \cdot)\|_{L^2}+
\|v'(t, \cdot)\|_{L^2}^{\f{1}{2}} 
\|v(t, \cdot)\|_{L^2}^{\f{3}{2}}).  
\end{eqnarray*}
We are almost  ready to derive an estimate on the desired norms, except that it is  hard to control $|v'(t,0)| $. Indeed, the easiest way out would be to use the norm $\|v'(t, \cdot)\|_{L^\infty}$, for which we unfortunately lack {\it a priori} estimates. Instead, we employ a simple averaging argument to reduce to the quantity $\|v'(t, \cdot)\|_{L^2}$, for which we do have  {\it a priori} bound, see \eqref{130}. 

Indeed, instead of finding the bound for $ \int v^2(t,x) |x| dx$, we might employ an estimate for $ \int v^2(t,x) |x-a| dx$, where $a\in (-1,1)$. Running through identical estimates, we arrive at 
\begin{eqnarray*}
	& & 	\p_t  \int v^2(t,x) |x-a| dx + \f{1}{4} \int v^2(t,x) dx \leq  C( |\dot{x}_0 (t)| \|v(t, \cdot)\|_{L^2}+|v(t,a)| |v'(t,a)| )+\\
	& & + C(\|v'(t, \cdot)\|_{L^2}^2+\|v'(t, \cdot)\|_{L^2} \|v(t, \cdot)\|_{L^2}+
	 \|v'(t, \cdot)\|_{L^2}^{\f{1}{2}} 
	\|v(t, \cdot)\|_{L^2}^{\f{3}{2}}).
\end{eqnarray*}
Integrating on the time interval $(0,t)$ brings about 
\begin{eqnarray*}
	& &  \int v^2(t,x) |x-a| dx  + \f{1}{4} \int_0^t \int v^2(s,x) dx ds \leq  \int v_0^2(x) |x-a| dx + \int_0^t |v(s,a)| |v'(s,a)| ds +\\
	& & + C(  \int_0^t |\dot{x}_0 (s)| \|v(s, \cdot)\|_{L^2}  + \|v'(s, \cdot)\|_{L^2}^2
	+\|v'(s, \cdot)\|_{L^2} \|v(s, \cdot)\|_{L^2}+ \|v'(s, \cdot)\|_{L^2}^{\f{1}{2}} 
	\|v(s, \cdot)\|_{L^2}^{\f{3}{2}} ) ds. 
\end{eqnarray*}
Taking an average of the previous estimate, i.e. integrating the previous estimate $\int_{-1}^1 \ldots   da$, and in doing so, recall  the formula 
$$
\int_{-1}^1 |x-a|   da=\left\{   \begin{array}{cc} 
	2|x| & |x|>1 \\ x^2+1 & |x|\leq 1
\end{array} \right. \in (2|x|, 2|x|+1).
$$
 This yields the estimate 
 \begin{eqnarray*}
 	& &  \int v^2(t,x) |x| dx  + \f{1}{2} \int_0^t \int v^2(s,x) dx ds \leq  
 	\int v_0^2(x) (2|x|+1) dx + \int_{-1}^1 \int_0^t |v(s,a)| |v'(s,a)| ) ds da+ \\
 	& & + C  \int_0^t [|\dot{x}_0 (s)| \|v(s, \cdot)\|_{L^2}  + \|v'(s, \cdot)\|_{L^2}^2+
 	\|v'(s, \cdot)\|_{L^2} \|v(s, \cdot)\|_{L^2}+ \|v'(s, \cdot)\|_{L^2}^{\f{1}{2}} 
 	\|v(s, \cdot)\|_{L^2}^{\f{3}{2}}] ds . 
 \end{eqnarray*}
Now, based on \eqref{130} and H\"older's inequality
\begin{eqnarray*}
& & 	 \int_0^t |\dot{x}_0 (s)| \|v(s, \cdot)\|_{L^2} ds \leq  \left( \int_0^t |\dot{x}_0 (s)|^2 ds \right)^{\f{1}{2}}  \left(\int_0^t \|v(s, \cdot)\|_{L^2}^2 ds\right)^{\f{1}{2}} \leq C \|v\|_{L^2_{t,x}}, \\
& &  \int_0^t \|v'(s, \cdot)\|_{L^2}^2 ds \leq C, \\
& & \int_0^t \|v'(s, \cdot)\|_{L^2} \|v(s, \cdot)\|_{L^2} ds\leq C \left(\int_0^t \|v'(s, \cdot)\|_{L^2}^2\right)^{\f{1}{2}} \left(\int_0^t \|v(s, \cdot)\|_{L^2}^2\right)^{\f{1}{2}} \leq C \|v\|_{L^2_{tx}},    \\
& &  \int_0^t  \|v'(s, \cdot)\|_{L^2}^{\f{1}{2}} 	\|v(s, \cdot)\|_{L^2}^{\f{3}{2}} ds   
 \leq \left(\int_0^t   \|v'(s, \cdot)\|_{L^2}^2 ds\right)^{\f{1}{4}} \left(\int_0^t   \|v(s, \cdot)\|_{L^2}^2 ds\right)^{\f{3}{4}}\leq C \|v\|_{L^2_{t,x}}^{\f{3}{2}},  \\
& &  \int_{-1}^1 \int_0^t |v(s,a)| |v'(s,a)|  da \leq  \left(\int_0^t   \|v'(s, \cdot)\|_{L^2}^2 ds\right)^{\f{1}{2}} \left(\int_0^t   \|v(s, \cdot)\|_{L^2}^2 ds\right)^{\f{1}{2}} \leq C \|v\|_{L^2_{t,x}}. 
\end{eqnarray*}
All in all, we conclude that
$$
 \int v^2(t,x) |x| dx+ \f{1}{2} \|v\|_{L^2_{t x}}^2 \leq C \|v\|_{L^2_{t,x}}^{\f{3}{2}}+C.
$$
By Young's inequality, we can hide $\|v\|_{L^2_{t,x}}^{\f{3}{2}}$ on the right-hand side behind $\f{1}{2} \|v\|_{L^2_{t,x}}^2$, which  leads us to the bound 
\begin{equation}
	\label{150} 
	\sup_{0<t<T} \int v^2(t,x) |x| dx+   \int_0^T \|v(t, \cdot)\|_{L^2_{x}}^2 dt \leq C.  
\end{equation}
for each $T>0$. This bound has two important consequences. 
First, using H\"older inequality, we have that for each $1<p<2$, 
\begin{eqnarray*}
	& &
\int |v(x)|^p dx = \int_{|x|<1} |v(x)|^p dx + \int_{|x|\geq 1} |v(x)|^p dx\leq  \\
& & \leq 
C \|v\|_{L^2}^p+ \left(\int_{|x|>1} v^2(x) |x| dx \right)^{\f{p}{2}} \left(\int_{|x|\geq 1} |x|^{-\f{p}{2-p}}dx \right)^{1-\f{p}{2}}\leq C_p  \left(\int |v|^2 (1+|x|) dx \right)^{\f{p}{2}},
\end{eqnarray*}
for which we have good bounds from \eqref{150}. It follows that 
 for any $p>1$, $\sup_{0<t< T} \int |v(x)|^p\leq C_p$. Note however that the constant $C_p$ blows up as $p\to 1+$, due to  $\lim_{p\to 1+} \left(\int_{|x|\geq 1} |x|^{-\f{p}{2-p}}dx \right)=+\infty$. 
 
  Second, by taking into account \eqref{129} and \eqref{150}, we have that for each $t>0$, 
$$
t \|v(t, \cdot)\|_{L^2}^2\leq \int_{0}^{t} \|v(s, \cdot)\|_{L^2}^2 ds \leq \int_0^{\infty} \|v(s, \cdot)\|_{L^2}^2 ds\leq C,
$$
whence it follows that $v$ obeys the bound $\|v(t, \cdot)\|_{L^2}\leq C t^{-\f{1}{2}}$. Interpolating this with the uniform bounds for $\|v(t, \cdot)\|_{L^p},  1<p<2$,  obtained just before, we have that for every $\de>0$, there is $C=C(\de)$, so that 
$$
\|v(t, \cdot)\|_{L^p}\leq C t^{-(1-\f{1}{p})+\de}, 1<p<2.
$$
\section{Asymptotic attraction to odd fronts for the generalized Burgers model} 
\label{sec:4} 
In this section, we present an extension of the results in \cite{BBHY}, when the fronts and the respective perturbations  are odd. Note that the fractional Burgers equation preserves odd initial data. 

Now that $\cl=-\sum_{j=1}^N a_j (-\p_{xx})^{\al_j} $ preserves the parity of the inputs, it is clear that if $u_0$ is odd, so is $v_0=u_0-\phi$. Clearly, the evolution preserves the odd solutions, whence $x\to v(t,x)$ is odd for all $t\geq 0$. 
Recall that by the construction in \cite{BBHY}, the translation function $x_0$ was selected (see formula $(20)$ there) via the formula  
$
\dot{x}_0(t)=\ga \dpr{\phi'}{v}=0,
$
since $\phi'$ is even, while $v$ is odd. It follows by Theorem \ref{theo:BBHY} that for this particular front and odd perturbations $u_0$,  $u(t,x)=\phi(x)+v(t,x)$, where 
\begin{equation}
	\label{49}
	\lim_{t\to \infty} \|v(t, \cdot)\|_{L^p}=0, 2< p\leq \infty.
\end{equation}
Recall that $v$ also satisfies the energy estimate \eqref{50}. 

We now aim to establish uniform estimates for the norm $\|v(t, \cdot)\|_{L^1(\rone)}$. To this end, take arbitrary large $T$ and assume $0\leq t<T$. Note that $v$ satisfies 
\begin{equation}
	\label{60} 
	v_t - v_{xx} - \cl v +\phi' v+\phi v'+ v v'=0.
\end{equation}
For each $p>1$, multiply \eqref{60}by $v |v|^{p-2}$ and integrate in $x$. 
We obtain, 
$$
\int v_t v |v|^{p-2} dx = \f{1}{p} \p_t \int |v|^p dx, 
$$
Integration by parts yields for each $p>1$, 
\begin{eqnarray*}
	& & \int  v v'v |v|^{p-2} dx =  \int v'	|v|^p dx = - \int v	(|v|^p)' dx =- p\int v v'v |v|^{p-2} dx,
\end{eqnarray*}
whence 
$$
\int  v v'v |v|^{p-2} dx=0. 
$$
By  \eqref{42}, we have the coercivity estimates
$$
-\int v_{xx}  v |v|^{p-2} dx \geq 0,   -\int (\cl v) v |v|^{p-2} dx =\sum_{j=1}^N a_j \int ((-\p_{xx})^{\al_j} v) v |v|^{p-2} dx \geq 0, 
$$
Next, 
$$
\int \phi v' v |v|^{p-2}  dx = \f{1}{p} \int  \phi (|v|^p)' dx = -  \f{1}{p} \int  \phi'  |v|^p dx. 
$$
All in all, we get 
\begin{equation}
	\label{60} 
	\p_t \int |v|^p dx +  (p-1) \int \phi' |v|^{p} dx \leq 0.
\end{equation}
Denoting $M_T:=\sup_{0<t<T} \|v(t, \cdot)\|_{L^\infty}$, and integrating the above expression in $(0,t)$, we obtain 
\begin{equation}
	\label{65}
	\int |v(t, x)|^p dx\leq \|v_0\|_{L^p}^p+ (p-1) \int_0^t \int_{\rone} |\phi'(x)| |v(t,x)|^{p} dx\leq  \|v_0\|_{L^p}^p+(p-1)  \|\phi'\|_{L^1}T M_T^p.
\end{equation}
Note that all the constants in \eqref{65} are  written out explicitly. This allows us to take a limit $p\to 1+$. We obtain 
\begin{equation}
	\label{70} 
	\|v(t, \cdot)\|_{L^1}\leq \|v_0\|_{L^1},
\end{equation}
and this is valid for all $T>0$, so \eqref{70} is a global estimate. 

This provides another important ingredient to the estimates \eqref{49}, which we have borrowed from Theorem \ref{theo:BBHY}. Fix $\eps>0$. By Plancherel's identity
$$
I(t):=\|v(t, \cdot)\|_{L^2}^2=\|v_{<\eps}(t, \cdot)\|_{L^2}^2+\|v_{>\eps}(t, \cdot)\|_{L^2}^2=:I_{<\eps}(t)+I_{>\eps}(t). 
$$
We now use \eqref{50} to conclude 
$$
\p_t(I_{<\eps}(t)+I_{>\eps}(t))\leq - C \|v_x(t)\|_{L^2}^2<- C \|\p_x v_{>\eps}(t)\|_{L^2}^2\leq -C \eps^2 I_{>\eps}(t).
$$
It follows that 
$$
I_{>\eps}'(t)+C\eps^2  I_{>\eps}(t)+  I_{<\eps}'(t)\leq 0.
$$
Multiplying by $e^{C \eps^2 t}$ and integrating yields 
\begin{equation}
	\label{79} 
	e^{C \eps^2 t} I_{>\eps}(t)  - I_{>\eps}(0) +\int_0^t e^{C \eps^2 s} I_{<\eps}'(s) ds\leq 0.
\end{equation}
However, integration by parts yields 
\begin{eqnarray*}
	-\int_0^t e^{C \eps^2 s} I_{<\eps}'(s) ds=-I_{<\eps} (s)e^{C \eps^2 s}|_0^t + C\eps^2 \int_0^t e^{C \eps^2 s} I_{<\eps}(s) ds\leq I_{<\eps} (0)+C\eps^2 \int_0^t e^{C \eps^2 s} I_{<\eps}(s) ds. 
\end{eqnarray*}
Plugging this back in \eqref{79} and multiplying by $e^{-C \eps^2 t}$, yields 
\begin{equation}
	\label{80} 
	I_{>\eps}(t)\leq I(0) e^{-C \eps^2 t}+ C\eps^2 \int_0^t e^{C \eps^2 (s-t)} I_{<\eps}(s) ds\leq I(0) e^{-C \eps^2 t}+\max_{0<s<t} I_{<\eps}(s).
\end{equation}
Now, we apply the  Bernstein's inequality for $A=(-\eps, \eps)$,  and \eqref{70} to obtain 
\begin{equation}
	\label{82} 
	I_{<\eps}(s)=
	\|v_{<\eps}(s, \cdot)\|_{L^2(\rone)}^2\leq C \eps \|v(s, \cdot)\|_{L^1(\rone)}^2 \leq C\eps \|v_0\|_{L^1(\rone)}^2. 
\end{equation}
Adding $I_{<\eps}(t)$ on both sides of \eqref{80} and using \eqref{82},  we obtain 
\begin{equation}
	\label{85}
	\|v(t, \cdot)\|_{L^2}^2=I(t)\leq e^{- C \eps^2 t} \|v_0\|_{L^2}^2+ C\eps \|v_0\|_{L^1}^2\leq 
	C_1 (e^{- C_1 \eps^2 t} +\eps)\|v_0\|_{L^2\cap L^1}^2,
\end{equation}
for some absolute  constant $C_1$, i.e.  independent on $\eps, t$. 

This inequality is true for every fixed $\eps>0$, and every $t$, but since clearly the best strategy is to optimize \eqref{85} in $\eps$ given $t$, we should proceed with some care. To this end, fix  large $t>>1$, and then select the unique 
$\eps<<1: e^{- C_1 \eps^2 t}=\eps$. Clearly, such choice of $\eps$ is possible, let us compute its asymptotic in terms of $t$. Since we need to solve, 
$$
\f{\eps^2}{\ln(\eps)}=-\f{1}{C_1 t},
$$
 its unique solution obeys the asymptotic 
\begin{equation}
	\label{700} 
	\eps\sim \sqrt{\f{\ln(t)}{t}}.
\end{equation}
It follows that 
\begin{equation}
	\label{100}
	\|v(t, \cdot)\|_{L^2}^2 \leq C_3 \f{\sqrt{\ln(t)}}{\sqrt{t}}   \|v_0\|_{L^2\cap L^1}^2 
\end{equation}
We conclude that $\|v(t, \cdot)\|_{L^2}\leq C \|v_0\|_{L^2\cap L^1} \left(\f{\ln(t)}{t}\right)^{\f{1}{4}}$, as stated.  Interpolating between the bounds \eqref{100} and \eqref{70}, the Gagliardo-Nirenberg's inequality now yields 
$$
\|v(t, \cdot)\|_{L^q}\leq C \|v_0\|_{L^2\cap L^1} 
\left(\f{\ln(t)}{t}\right)^{\f{1-\f{1}{q}}{2}}, 1<q<2, 
$$
for intermediate values of $q$. 
While one can similarly interpolate with the bound \eqref{49}, this will not be a particularly good decay, especially near $q=\infty$. 

Instead, we use an alternative estimate based again on \eqref{50}. Indeed from  \eqref{130}  and \eqref{50} , it is clear that for every $T>0$ 
$$
\int_T^\infty \|v'(s, \cdot)\|_{L^2}^2 ds  \leq \|v(T, \cdot)\|_{L^2}^2\lesssim \left(\f{\ln(T)}{T}\right)^{\f{1}{2}},
$$
where the last estimate follows from the $L^2$ bounds \eqref{100}, 
Also, for every $T>>1$, we have that there exists $T_0\in (T, 2T)$, so that 
$$
T  \|v'(T_0, \cdot)\|_{L^2}^2 \leq \int_T^{2T}  \|v'(s, \cdot)\|_{L^2}^2 ds<\int_T^\infty \|v'(s, \cdot)\|_{L^2}^2 ds\lesssim \left(\f{\ln(T)}{T}\right)^{\f{1}{2}}. 
$$
In other words, 
\begin{equation}
	\label{110} 
	\|v'(T_0, \cdot)\|_{L^2}^2 \leq C \f{\sqrt{\ln(T_0)}}{T_0^{\f{3}{2}}}, \int_{T_0}^\infty \|v'(s, \cdot)\|_{L^2}^2 ds \leq  C \left(\f{\ln(T_0)}{T_0}\right)^{\f{1}{2}},
\end{equation}

For large $T_0$, denote $\de:=\left(\f{\ln(T_0)}{T_0}\right)^{\f{1}{2}}>>C \f{\sqrt{\ln(T_0)}}{T_0^{\f{3}{2}}}$. A this point, we refer to the argument in \cite{BBHY}, which establishes that starting with \eqref{110}, we conclude 
\begin{equation}
	\label{115} 
	\sup_{t>T_0} \|v'(t)\|_{L^2}\leq C \de^{\f{2}{3}}\sim \left(\f{\ln(T_0)}{T_0}\right)^{\f{1}{3}}
\end{equation}
Since this can be done for at least one point $T_0\in (2^k, 2^{k+1}), k>>1$, we obtain the bound 
$$
\|v'(t)\|_{L^2}\leq C \left(\f{\ln(t)}{t}\right)^{\f{1}{3}}
$$
By Sobolev embedding, followed by the Gagliardo-Nirenberg's inequality  we  conclude 
$$
\|v(t, \cdot)\|_{L^q}\leq \|v(t, \cdot)\|_{\dot{H}^{\f{1}{2}-\f{1}{q}}}\leq \|v(t, \cdot)\|_{L^2}^{\f{1}{2}+\f{1}{q}}\|v'(t, \cdot)\|_{L^2}^{\f{1}{2}-\f{1}{q}}\leq  C \left(\f{\ln(t)}{t}\right)^{\f{7}{24}-\f{1}{12 q}}.
$$



\begin{thebibliography}{99}
	
	\bibitem{BBHY} B. Barker, J., Bronski, V. Hur, Z. Yang, {\emph Asymptotic stability of sharp fronts. I. One bound state implies stability}, (2022), available  \href{https://arxiv.org/pdf/2112.04700.pdf}{here}
	
\bibitem{BBHY2} B. Barker, J., Bronski, V. Hur, Z. Yang, {\emph	Asymptotic stability of sharp fronts: Analysis and rigorous computation}, available  \href{https://arxiv.org/pdf/2112.04700}{here}
	
\bibitem{BS1} 	J. Bona, M. Schonbek, {\emph Traveling wave solutions to Korteweg-de Vries-Burgers equation.} {\em  Proc Roy Soc
Edinburgh}, {\bf 101A}, (1985), p. 207--226.
	
\bibitem{BVS} 	J. Bona, S. Rajopadhye, M. E. Schonbek, {\emph Models for propagation of bores. I. Two-dimensional theory.} {\em Diff. Int. Eq.}, {\bf 7}, No. 3,4 (1994), p. 699--734. 

\bibitem{CL} D. Chamorro,  P. G. Lemari\'e-Rieusset, {\emph Quasi-geostrophic equations, nonlinear Bernstein inequalities and $\al$-stable processes.} {\em  Rev. Mat. Iberoam.} {\bf  28},  (2012), no. 4,  p. 1109--1122.

\bibitem{JM} A. Jeffrey  M.. Mohamad   {\emph Exact solutions to the KdV-Burgers equation},  {\em Wave Motion}, {\bf 14}, (1991), p. 369--375. 

\bibitem{NS1} P. I. Naumkin,  I. A. Shishmarev {\emph A problem on the decay of step-like data for the KortewegdeVries-Burgers equation}, {\em  Funktsional. Anal. i Prilozhen}, {\bf 25}, (1):(1991), p. 21--32.
\bibitem{NS2} P. I. Naumkin, I. A. Shishmarev {\emph On the decay of step-like data for the Korteweg-deVries-Burgers equation}, {\em  Funktsional. Anal. i Prilozhen},  {\bf 26}, (2), (1992), p. 88--93. 

\bibitem{Pego} R. L. Pego. {\emph Remarks on the stability of shock profiles for conservation laws with dissipation.} {\em Trans. Amer. Math. Soc.}, {\bf  291} (1), (1985), p. 353--361. 
 
 
 \bibitem{Stein} E. M. Stein,  {\emph Topics in harmonic analysis related to the Littlewood–Paley theory.} {\em 
Annals of Mathematics Studies}, {\bf  63}, 
Princeton University Press, Princeton, NJ; University of Tokyo Press, Tokyo, 1970.
	
\bibitem{SG} 	C.  Su, C.  Gardner, {\emph Korteweg-de Vries equation and generalizations. III. Derivation
of the Korteweg-de Vries equation and Burgers equation.} {\em J. Math. Phys.}  
{\bf 10}, (3)  (1969), p. 536--539,
	
\bibitem{ZH}	K. Zumbrun, P. Howard, {\emph Pointwise semigroup methods and stability of viscous shockwaves}, {\em  Indiana Univ. Math. J.}, {\bf 47}, (3), (1998) p. 741--871. 
	
	
	
	
%
%
%
%
%
%
%
%
%
%
%
%
%
%
%
%
%
%
%
%
%
%
%
%
%
%
%
%
%
%
%
%
%
%
%
%
%
%
%
%
%
%
%
%
%
%
%
%
%
%
%
%
%
%
%
%
%
%

\end{thebibliography}
\end{document}